\documentclass{amsart}

\usepackage{amsmath,amsthm,amsfonts,latexsym,euscript}

\title [Structure and isomorphism classification of $A_u(Q)$ and $B_u(Q)$]
{\bf Structure and Isomorphism Classification of
Compact Quantum Groups $A_u(Q)$ and $B_u(Q)$}
\author[Shuzhou Wang]
{\bf Shuzhou Wang}
\address{Department of Mathematics, University of Georgia,  
Athens, GA 30602 
\newline \indent
Fax: 706-542-2573  
}
\email{szwang@@math.uga.edu} 

\date{}

\newtheorem{DF}{Definition}[section]
\newtheorem{LM}[DF]{Lemma}
\newtheorem{PROP}[DF]{Proposition}
\newtheorem{THM}[DF]{Theorem}
\newtheorem{COR}[DF]{Corollary}
\newtheorem{RMK}[DF]{Remark}
\newtheorem{RMKS}[DF]{Remarks}
\newtheorem{PROB}[DF]{Problem}

\newcommand{\bgdf}{\begin{DF}}
\newcommand{\nddf}{\end{DF}}
\newcommand{\bglm}{\begin{LM}}
\newcommand{\ndlm}{\end{LM}}
\newcommand{\bgprop}{\begin{PROP}}
\newcommand{\ndprop}{\end{PROP}}
\newcommand{\bgthm}{\begin{THM}}
\newcommand{\ndthm}{\end{THM}}
\newcommand{\bgcor}{\begin{COR}}
\newcommand{\ndcor}{\end{COR}}
\newcommand{\bgrmk}{\begin{RMK}}
\newcommand{\ndrmk}{\end{RMK}}
\newcommand{\bgrmks}{\begin{RMKS}}
\newcommand{\ndrmks}{\end{RMKS}}
\newcommand{\bgprob}{\begin{PROB}}
\newcommand{\ndprob}{\end{PROB}}
\newcommand{\bgeq}{\begin{eqnarray}}
\newcommand{\ndeq}{\end{eqnarray}}
\newcommand{\bgeqq}{\begin{eqnarray*}}
\newcommand{\ndeqq}{\end{eqnarray*}}

\newcommand{\pf}{{\em Proof. }} 
\newcommand{\QED}{\hfill Q.E.D.} 
\newcommand{\vv}{\vspace{4mm}\\}
\newcommand{\vvv}{\vspace{6mm}\\}

\def\CC{{\mathbb C}}

\numberwithin{equation}{section}

\begin{document}

\begin{abstract}
We classify the compact quantum groups
$A_u(Q)$ (resp. $B_u(Q)$) up to isomorphism 
when $Q>0$ (resp. when $Q \bar{Q} \in {\mathbb R} I_n$).  
We show that the general $A_u(Q)$'s and $B_u(Q)$'s for 
arbitrary $Q$ have explicit decompositions into 
free products of these special $A_u(Q)$'s and $B_u(Q)$'s.  
\end{abstract}

\maketitle

\noindent 
\vvv
{\bf Introduction}
\vv
Recall \cite{Wor5,Wor9} that a {\em compact matrix quantum group} is a pair 
$G=(A,u)$ of a unital $C^*$-algebra $A$ and matrix system
$u$ of generators $u_{ij}$   ($i,j=1, \cdots, n$)
that satisfies the following two axioms:

  (1) There is a unital $C^*$-homomorphism $\Phi: \; A
\longrightarrow A \otimes A$ such that
$\Phi(u_{ij})
=\sum ^{n}_{k=1} u_{ik} \otimes u_{kj}$
for each  $i,j$;

  (2) The matrices $u=(u_{ij})$ and
  $u^t$ are invertible in $M_{n}(\CC) \otimes A$.

In \cite{W5,W5',W1}, we constructed
for each $Q \in GL(n, {\mathbb C})$ two families
$A_u(Q)$ and $B_u(Q)$ of compact matrix quantum groups
in the sense of Woronowicz \cite{Wor5}.
The compact
quantum group $A_u(Q)$ and $B_u(Q)$ are is defined in terms of generators
$u_{ij}$ ($i,j =1,  \cdots  n$), and relations: 
\bgeqq
& & A_u(Q): \; \; \;
    u^* u = I_n = u u^*, \; \; \;
    u^t Q {\bar u} Q^{-1} = I_n = Q {\bar u } Q^{-1} u^t; \\
& & B_u(Q): \; \; \;
    u^* u = I_n = u u^*, \; \; \;
    u^t Q u Q^{-1} = I_n = Q u Q^{-1} u^t,
\ndeqq
where $u = (u_{ij})$. The $A_u(Q)$'s are universal in the sense 
that every compact matrix quantum group
is a quantum subgroup of $A_u(Q)$ for some $Q>0$. 
Similarly, the $B_u(Q)$'s are universal in the sense 
that every compact matrix quantum group with self conjugate
     fundamental representation 
is a quantum subgroup of $B_u(Q)$ for some $Q$. 
The subscript $u$ denotes ``universality''. 
For $Q>0$ (resp. $Q$ with $Q \bar{Q} \in {\mathbb R} I_n$), 
Banica determined in \cite{Banica2} (resp. \cite{Banica1}) the fusion rings of 
the irreducible representations of the quantum group $A_u(Q)$ 
(resp. $B_u(Q)$). Note that he used 
$A_u(F)$  (resp. $A_o(F)$) to denote $A_u(Q)$ (resp. $B_u(Q)$) 
where $Q = F^* F$ (resp. $Q = F^*$). 
Although these quantum groups are of a different nature 
from those well known ones obtained from ordinary
Lie groups by deformation quantization \cite{Soib1,VS2,Wor4,Wor6}, they share
many properties of the ordinary Lie groups. Not only they are
fundamental objects in the framework of Woronowicz on compact matrix
quantum groups \cite{Wor5}, they are also useful objects in the study
of intrinsic quantum group symmetries, such as ergodic quantum group
symmetries on operator algebras and quantum automorphism groups
of noncommutative spaces (cf. \cite{W14,W15}).
Although much is known about $A_u(Q)$ and $B_u(Q)$, some basic
questions are still left unanswered \cite{W7}. For instance, for
different $Q$'s, how do the $A_u(Q)$'s (resp. $B_u(Q)$'s)
differ from each other? Do some of these quantum groups constitute building
blocks for the two families of quantum groups in an appropriate sense?

The purpose of this paper is to answer these questions.
Recall that the fundamental representation $u$ of $A_u(Q)$ (resp. $B_u(Q)$) is
irreducible if and only if $Q$ is positive (resp.
$Q \bar{Q} \in {\mathbb R} I_n$), see \cite{W14} (resp. \cite{Banica1}).
When these conditions are satisfied, we classify $A_u(Q)$ and
$B_u(Q)$ up to isomorphism, and show that they are not free
products or tensor products or crossed products (see \cite{W1,W2}
for these constructions). 
We show that the general $A_u(Q)$'s and $B_u(Q)$'s for 
arbitrary $Q$ are free products of 
these special $A_u(Q)$'s and $B_u(Q)$'s, 
and we give their explicit decompositions
in terms of free products. 

In the following, the word morphisms means morphisms between compact quantum
groups (cf. \cite{W1}).
\vvv
{\bf 1. The quantum groups $A_u(Q)$ for positive $Q$}
\vv
Let $Q \in GL(n, {\mathbb C})$. Then $A_u(Q)=A_u(cQ)$ for any
nonzero number $c$. For a positive matrix $Q$, we can
normalize it so that $Tr(Q) = Tr(Q^{-1})$.
\vv
THEOREM 1.
{\em Let $Q \in GL(n, {\mathbb C})$ and $Q' \in GL(n', {\mathbb C})$
be positive matrices normalized as above with eigen values
$q_1 \geq q_2 \geq \cdots \geq q_n$ and
$q'_1 \geq q'_2 \geq \cdots \geq q'_{n'}$ respectively. Then,}

(1). {\em 
$A_u(Q)$ is isomorphic to $A_u(Q')$ if and only if,
(i)
$n = n'$, and
(ii)
$  (q_1, q_2, \cdots, q_n)$
$= (q'_1, q'_2, \cdots, q'_n) $ or
$(q^{-1}_n, q^{-1}_{n-1}, \cdots, q^{-1}_1)= (q'_1, q'_2, \cdots, q'_n)$.}

(2). {\em 
$A_u(Q)$ is not a free product. That is, if $A_u(Q)= A \ast B$ is a free
product of Woronowicz $C^*$-algebras $A$ and $B$, then either
$A= A_u(Q)$ or $B= A_u(Q)$.}
\vv
\pf Clearly, we may assume $n, n' \geq 2$.

(1).
Since $A_u(Q)$ is similar to
$A_u(VQV^{-1})$ for $V \in U(n)$ (cf. \cite{W5}), we may assume up to
isomorphism that $Q$ and $Q'$ are in diagonal form, say,
$$Q = diag(q_1, q_2, \cdots, q_n), \hspace{1cm}
Q' = diag(q'_1, q'_2, \cdots, q'_{n'}).$$
We claim that every non-trivial
irreducible representation of the quantum group $A_u(Q)$ other than $u$ and
$\bar{u}$ has a dimension greater than $n$.

In \cite{Banica2}, the irreducible representations of
the quantum group $A_u(Q)$ are parameterized by
the free monoid ${\mathbb N} \ast {\mathbb N}$ with generators
$\alpha$ and $\beta$ and anti-multiplicative involution $\bar{\alpha}= \beta$.
(the neutral element is $e$ with $\bar{e} = e$).
The classes of $u$ and $\bar{u}$ are $\alpha$
and $\beta$ respectively.
Let $d_x$ be the dimension of irreducible representations in the class
 $x \in {\mathbb N} \ast {\mathbb N}$.
By Theorem 1 of \cite{Banica2}, we have the following dimension 
formula: 
$$d_x d_y = \sum_{x = ag, \bar{g} b = y}  d_{ab}. $$
Hence $d_{x \alpha^{k + l} y} = d_{x \alpha^k} d_{\alpha^l y}$
for $k, l \geq 1$. This identity prevails
when we change $\alpha$ to $\beta$.
From these we see that apart from the trivial class $e$,
the classes with smaller
dimensions are those words $x$ in which the powers of $\alpha$ and $\beta$
are equal to 1. Moreover, we infer from \cite{Banica2} that
for any word $x$, $d_x$ does not change when we exchange
$\alpha$ and $\beta$ in $x$. Hence the minimal dimension
is among $d_{\alpha}, d_{\alpha \beta}, d_{\alpha \beta \alpha}, \cdots,$ 
so we now concentrate on these numbers.  
Let $f(1), f(2), f(3), \cdots$ be this sequence and let $f(0)=1 =d_e$. Then 
applying the above dimension formula to 
$d_{\alpha \beta... \alpha \beta} d_\alpha $
and 
$d_{\alpha \beta... \alpha \beta \alpha} d_\beta, $
 we get   
$f(k+1) = n f(k)- f(k-1)$, $k \geq 1$,  
noting  that $f(1)= d_\alpha = d_\beta = n$. Since $n \geq 2$,
we have
\bgeqq
f(k+1) - f(k) &  =  & (n-1) f(k) - f(k-1) \geq
f(k) - f(k-1) \geq \cdots \\
              & \geq & f(1) - f(0) >0.
\ndeqq
Hence $d_\alpha = d_\beta = n < d_x$ for
$x \neq \alpha, \beta, e$. This proves our claim.

We introduced in \cite{W6}
$F$-matrices for classes of irreducible representations of a
compact quantum group, based on Woronowicz \cite{Wor5}.
If $v$ is an irreducible representation
with $F$-matrix $F_v$ in the sense of \cite{Wor5}, then the $F$-matrix
$F_{[v]}$ for the class $[v]$ of $v$ is the diagonal matrix
with eigen values of $F_v$ arranged in decreasing
order on the diagonal. This is an invariant of the class $[v]$.
We have $F_{[\bar{v}]} = diag(\lambda^{-1}_m, \lambda^{-1}_{m-1},
\cdots, \lambda^{-1}_1)$ if $F_{[v]} = diag(\lambda_1, \lambda_2, \cdots,
\lambda_m)$ (Warning: $F_{\bar{v}} = (F_v^t)^{-1}$).

Now for the quantum group $A_u(Q)$,
$F_{[u]} = F_u = Q^t = Q$ (cf. Remark 1.5.(3) of \cite{W5}).
Therefore if the quantum groups $A_u(Q)$ and $A_u(Q')$
are isomorphic to each other, the $F$-matrices for the {\em classes}
of the irreducible representations $u$ and $\bar{u}$ (of minimal
dimension $n$) of $A_u(Q)$ correspond to those of $A_u(Q')$ (cf. Lemma 4.2 of
\cite{W6}). Whence we have conditions (i) and (ii) in the theorem.
Conversely, assume conditions (i) and (ii) are satisfied. 
If $Q = Q'$, there is nothing to prove. 
So we assume that 
$(q^{-1}_n, q^{-1}_{n-1}, \cdots, q^{-1}_1)= (q'_1, q'_2, \cdots, q'_n)$.
As $F_{[\bar{u}]} = diag (q^{-1}_n, q^{-1}_{n-1}, \cdots, q^{-1}_1)$,
we can choose a unitary representation $v$ in the class $[\bar{u}]$
of $\bar{u}$ (non-unitary in general)
so that $F_v = F_{[\bar{u}]}$. Then the entries of
$v$ generate the same algebra $A_u(Q)$ as those of $u$ and
they satisfy the relations for $A_u(Q')$. It is now clear that 
$A_u(Q)$ and $A_u(Q')$ are isomorphic to each other,

(2). Suppose $A_u(Q) = A \ast B$. By the classification of irreducible
representations of the quantum $A \ast B$ in \cite{W1},
the representation $u$ is a tensor product of non-trivial
irreducible representations of the quantum groups $A$ and $B$. Also by
\cite{W1}, each representation in this tensor product is also an
irreducible representation of $A_u(Q)$.
From the claim in the proof of (1) above,
we deduce that $A_u(Q)$ has no irreducible representation of dimension 1
other than the trivial one.
Therefore there is only one term in the tensor product.
That is, $u$ is a representation of the quantum group $A$ or $B$.
Whence $A_u(Q)=A$ or $A_u(Q)=B$.
\QED
\vv
{\em Remarks.}
(1). As a corollary of the proof above,
we have the following rigidity result for
$A_u(Q)$ (a similar result holds for $B_u(Q)$ in the next section).
 Let $Q, Q' \in GL(n, {\mathbb C})$ be positive, normalized as above.
If $A_u(Q')$ is a quantum subgroup of $A_u(Q)$ given by a surjection
$\pi: A_u(Q) \rightarrow A_u(Q')$, then $Q' = VQV^{-1}$
or $Q' = V(Q^t)^{-1}V^{-1}$ for some $V \in U(n)$ and hence $\pi$ is
an isomorphism. To see this, first
$S u' S^{-1} = \pi(u)$ or
    $S {Q'}^{1/2} \overline{u'} {Q'}^{-1/2} S^{-1} = \pi(u)$
     for some $V \in U(n)$, which satisfy the relations for $A_u(Q)$.
The assertion follows from the irreducibility of the representation $u'$.

(2).
By the same method, one can also prove that $A_u(Q)$ is not
a tensor product, nor a crossed product (cf. \cite{W2}). This remark
also applies to $B_u(Q)$ blow.

(3).
Although $A_u(Q)$ is a universal analog of $U(n)$ for compact quantum groups
\cite{W1,W5',W5}, the proof in the above shows that $A_u(Q)$
has no nontrivial irreducible representations of dimension 1 when $n \geq 2$,
a property in sharp contrast to $U(n)$. Further study of
the irreducible representations of the quantum group $A_u(Q)$
gives evidences that $A_u(Q)$ may be a {\em simple} compact
quantum group in an appropriate sense (work in progress).
\vvv
{\bf 2. The quantum groups $B_u(Q)$ with
$Q \bar{Q} \in {\mathbb R} I_n$}
\vv
The quantum group $B_u(Q)$ has only
one irreducible representation of minimal dimension 
among the non-trivial ones (cf. \cite{Banica1}).
If $B_u(Q)$ is isomorphic to $B_u(Q')$, then the fundamental
representation of $B_u(Q)$ corresponds to that of $B_u(Q')$ under the
isomorphism. Using the irreducibility of the fundamental representations
along with the defining relations for $B_u(Q)$ and $B_u(Q')$
we immediately obtain part (1) of the following theorem
(we need not consider the $F$-matrices for this).
The proof of part (2) of the theorem
is similar to the proof of part (2) of Theorem 1 above and will be omitted.
\vv
THEOREM 2.
{\em Let $Q \in GL(n, {\mathbb C})$ and $Q' \in GL(n', {\mathbb C})$
be matrices such that $Q \bar{Q} \in  {\mathbb R} I_n$ and
$Q' \overline{Q'} \in  {\mathbb R} I_{n'}$. Then,}

(1). {\em $B_u(Q)$ is isomorphic to $B_u(Q')$ if and only if,
(i)
$n = n'$, and
(ii)
there exist $S \in U(n)$ and $c \in {\mathbb C}^*$ such that
$Q = z S^t Q' S$.}

(2). {\em 
$B_u(Q)$ is not a free product. That is, if $B_u(Q)= A \ast B$ is a free
product of Woronowicz $C^*$-algebras $A$ and $B$, then either
$A= B_u(Q)$ or $B= B_u(Q)$.}
\vv
{\bf Explicit parametrization of the isomorphism classes of $B_u(Q)$}.
Contrary to Theorem 1, Theorem 2 above does not give an explicit
parametrization of the isomorphism classes of the $B_u(Q)$'s. 
Assume the normalization $Tr(Q Q^*)$ $= Tr((Q Q^* )^{-1})$ 
(or, equivalently, $Tr(\bar{Q} Q^t)$ $= Tr((\bar{Q} Q^t )^{-1})$). 
Let $\lambda_1 \geq \lambda_2 \geq \cdots \geq \lambda_n$
be the eigen values of $Q Q^*$ (therefore of $\bar{Q} Q^t$). 
Using the defining relations of $B_u(Q)$, 
we have that the antipode $\kappa$ satisfies  
$$\kappa^2(u) = \kappa(Q^{-1} u^t Q) = Q^{-1} Q^t u (Q^{-1})^t Q,$$
where $u$ is the fundamental representation of
$B_u(Q)$. Hence from the assumption that 
$Q \bar{Q} = cI_n$ we get  
$F_u = \bar{Q} Q^t$  (also cf. Remark 1.5.(3) of \cite{W5}).
Since $u = Q^* \bar{u} {Q^*}^{-1}$,
one has $F_{[u]}=F_{[\bar{u}]}$ and hence
$$F_{[u]} = (\lambda_1, \lambda_2, \cdots, \lambda_n)
=(\lambda_n^{-1}, \lambda_{n-1}^{-1}, \cdots, \lambda_1^{-1})=F_{[\bar{u}]}.$$
Let  $Q^* = U |Q^*|= U \sqrt{Q Q^*}$ 
be the polar decomposition. We have that 
$$c I_n = Q \bar{Q} = |Q^*| U^*  |Q^*|^t U^t.$$
Taking determinants of both sides of 
$$c I_n |Q^*|^{-1} =   U^*  |Q^*|^t U^t,$$ 
we get $c^n = 1$. Note that $c$ is real 
($c I_n= Q \bar{Q} =  \bar{Q} Q = \overline{Q \bar{Q}}$), 
we have $c =\pm 1$. 
Taking determinant on both sides of $Q \bar{Q} = cI_n$, 
we see that $c = \pm 1$ for $n$ even and $c=1$ for $n$ odd. 

Conversely, we claim that $ Q \bar{Q} = c I_n $ 
with $c =\pm 1$ for $n$ even and $c=1$ for 
$n$ odd implies that $Tr(Q Q^*)$ $= Tr((Q Q^* )^{-1})$. 
To see this, let $r>0$ be such that 
$Tr(Q_1 Q_1^*)$ $= Tr((Q_1 Q_1^* )^{-1}),$ 
where $Q_1 = r Q$. Then the above analysis applied to 
$B_u(Q_1)$ yields $ Q_1 \bar{Q_1} = c I_n $, 
i.e. $r^2 Q \bar{Q} = c I_n $. Hence $r =1$. 

Let us summarize this discussion in the following
\vv
PROPOSITION 1.
{\em Let $Q \in GL(n, {\mathbb C})$
be a matrix  such that $Q \bar{Q} \in  {\mathbb R} I_n$. Then
the condition $Tr(Q Q^*)$ $= Tr((Q Q^* )^{-1})$ is equivalent to
$Q \bar{Q} = \pm I_n$ for $n$ even and
$Q \bar{Q} =  I_n$ for $n$ odd. }
\vv
It would be interesting to find an elementary proof of the above fact 
(i.e. a proof without using quantum group theory).

Since $B_u(rQ) =B_u(Q)$ for any non-zero $r$, we 
will assume the normalization $Q \bar{Q} = \pm I_n$ below to find  
a parametrization of the equivalent classes of $B_u(Q)$. 
We note that if both $Q$ and $Q'$ are so normalized and that 
$Q' = z S^t Q S$ for some non-zero complex number $z$, then a straightforward 
computation shows that $z$ has modulus $1$. So  
we can restrict the $z$  in Theorem 2.(1) 
to complex numbers of modulus $1$. 

We will need the following easy lemma. 
\vv
LEMMA 1. {\em 
Let $Q' = z S^t Q S$, where $z$ is a number of modulus $1$, $S$ is a 
unitary scalar matrix of the same size as $Q$. Let 
$Q= U |Q|$ and $Q'= U' |Q'|$ be the polar decompositions of $Q$ and 
$Q'$. Then $|Q'| = S^{-1} |Q| S $ and $U' = z S^t U S $.
}
\vv
\pf
The identity $|Q'| = S^{-1} |Q| S$ follows from ${Q'}^*Q' = S^{-1} Q^* Q S$. 

Now, on the one hand we have 
$$Q' = U' |Q'| = U' S^{-1} |Q| S, $$
on the other hand we have 
$$Q' = z S^t Q S = z S^t U |Q| S. $$
Hence $U' S^{-1} |Q| S = z S^t U |Q| S ,$ and $U' =z S^t U S$.
\QED
\vv
From the analysis preceding Proposition 1, we can assume that  
the eigen values of $|Q|$ has the form 
$$\mu_1 \geq \mu_2 \geq \cdots \geq \mu_k 
 \geq \mu_k^{-1} \geq \cdots \geq \mu_2^{-1} \geq \mu_1^{-1}
$$
or  
$$\mu_1 \geq \mu_2 \geq \cdots \geq \mu_k 
\geq 1 \geq \mu_k^{-1} \geq \cdots \geq \mu_2^{-1} \geq \mu_1^{-1}
$$
according to $n = 2k$ or $n =2k +1$. 
In virtue of of Lemma 1 and Theorem 2.(1), we can assume that 
$$|Q| =diag(\mu_1, \mu_2, \cdots, \mu_k, \mu_k^{-1}, 
\cdots, \mu_2^{-1}, \mu_1^{-1})$$
or
$$|Q| =diag(\mu_1, \mu_2, \cdots, \mu_k, 1, 
\mu_k^{-1}, \cdots, \mu_2^{-1}, \mu_1^{-1})$$
is in diagonal form. 
So $Q = U diag(\mu_1, \mu_2, \cdots,  \mu_2^{-1}, \mu_1^{-1})$.  
Then from the normalization $Q \bar{Q} = \pm I_n$ and Lemma 1, 
we immediately obtain the following more explicit form of Theorem 2.(1):
\vv
THEOREM 2$^\prime$. {\em
The isomorphism classes of $B_u(Q)$ are given by
\begin{center}
$(U, (\mu_1, \mu_2, \cdots , \mu_k))$, \hspace{5mm}
$\mu_1 \geq \mu_2 \geq \cdots \geq \mu_k \geq 1$,
\end{center}
where $U \in U(n)$ is a solution of the equation
\begin{center}
  $U diag(\mu_1, \mu_2, \cdots , \mu_k, 
\mu_k^{-1}, \cdots, \mu_2^{-1}, \mu_1^{-1})$
  $=
c \cdot diag(\mu_1^{-1}, \mu_2^{-1}, \cdots, 
\mu_k^{-1}, \mu_k, \cdots , \mu_2, \mu_1) U^t,$
\end{center}
 or the equation
\begin{center}
$U diag(\mu_1, \mu_2, \cdots, \mu_k, 1, 
\mu_k^{-1}, \cdots, \mu_2^{-1}, \mu_1^{-1})$
$=
diag(\mu_1^{-1}, \mu_2^{-1}, \cdots, 
\mu_k^{-1}, 1, \mu_k, \cdots, \mu_2, \mu_1)
U^t$
\end{center}
according to $n = 2k$ ($c = 1$ or $-1$) or $n = 2k+1$.  
The pairs $(U, (\mu_1, \mu_2, \cdots , \mu_k))$ and
$(U', (\mu'_1, \mu'_2, \cdots , \mu'_k))$  
represent the same class if and only if
$(\mu_1, \mu_2, \cdots , \mu_k) = (\mu'_1, \mu'_2, \cdots , \mu'_k)$ and
$U' = z S^t U S$
for a $z \in {\mathbb T}$ and a stabilizing unitary $S$:
\begin{center}
$S
 diag(\mu_1, \mu_2, \cdots, \mu_2^{-1}, \mu_1^{-1}) S^{-1} =
 diag(\mu_1, \mu_2, \cdots, \mu_2^{-1}, \mu_1^{-1})$.
\end{center}
}
We can now easily recover the classification
in \cite{W6} for $SU_q(2)$ from the above. 
First we note that  
$C(SU_q(2)) = B_u(Q)$ with normalized 
$Q =
\left[
\begin{array}{cc}
                  0 & - s \sqrt{|q|}^{-1}    \\  
         \sqrt{|q|} & 0
\end{array}
\right]$,
where $s = q^{-1} |q|$ 
(cf. Sect. 5 of \cite{Banica1} or \cite{W5'}). 
In this case $Q = U |Q|$ with 
$$|Q| = 
\left[
\begin{array}{cc}
             \sqrt{|q|}^{-1} & 0    \\  
                           0 & \sqrt{|q|}
\end{array}
\right], \; \; \; \; \; 
U = 
\left[
\begin{array}{cc}
             0 & -s    \\  
             1 & 0
\end{array}
\right] .$$
 So the parametrization for 
$C(SU_q(2))$ in terms of Theorem 2$^\prime$ is 
$(U, \sqrt{|q|}^{-1})$, which is equivalent to 
saying that they are non-isomorphic to each other for
$q \in [-1,1] \backslash \{ 0 \}$ 
(cf.  Theorem 3.1 of \cite{W6}). 
 
Our work in progress shows that the quantum groups $B_u(Q)$ are
{\em simple} when $Q \bar{Q} \in {\mathbb R} I_n$.
\newpage
\noindent
{\bf 3. The quantum groups $A_u(Q)$ and $B_u(Q)$ for arbitrary $Q$}
\vv
In the following, $h$ will denote the Haar measure of the
ambient quantum group.
Note that $A_u(Q)=C({\mathbb T})$ and 
$B_u(Q) = C^*({\mathbb Z}/2{\mathbb Z})$ for
$Q \in GL(1, {\mathbb C})$.
\vv
THEOREM 3.
{\em Let $Q \in GL(n, {\mathbb C})$. Let
$u = S diag (m_1 w_1, m_2 w_2, \cdots, m_k w_k) S^{-1}$
be a decomposition of $u$ of $A_u(Q)$ into unitary
isotypical components $m_j w_j$ ($j=1, \cdots ,k$) 
for some $S \in U(n)$. Let $Q_j$ be the positive matrix
$h (w_j^t \bar{w}_j)$. 
Then
\begin{center}
$A_u(Q) \cong A_u(Q_1) \ast A_u(Q_2) \ast \cdots \ast A_u(Q_k).$
\end{center}
}
\noindent
\pf
Let $E = S^t Q \bar{S}$. 
 Then the second set of 
relations for $A_u(Q)$ becomes 
$$ w^t E \bar{w} E^{-1} = I_n = E \bar{w} E^{-1} w^t, \hspace{.51cm}
{\mbox i.e.} \hspace{.51cm}
w^t E \bar{w} = E, 
$$
where 
$w=diag(v_1, v_2, \cdots, v_k)$,  
 $v_j = m_j w_j$, $j = 1, \cdots, k$.
Block decompose $E$ according to
$w$, say, $E=(E_{ij})_{i,j=1}^k$.
 Then the above set of relations 
becomes $v_i^t E_{ij} \bar{v}_j = E_{ij}$ ($i,j=1, \cdots ,k$).
 By Lemma 1.2 of \cite{W5}, $(v^t_i)^{-1} = \tilde{Q}_i \bar{v}_i
\tilde{Q}_i^{-1}$,
where $\tilde{Q}_i = diag(Q_i, \cdots , Q_i)$ ($m_i$ copies). Hence
$\tilde{Q}_i^{-1} E_{ij} \bar{v}_j = \bar{v}_i \tilde{Q}_i^{-1} E_{ij} $.
Since the $w_i$'s are mutually inequivalent irreducible representations,
we deduce that $E_{ij}=0$ for $i \neq j$, and that $E_{jj}$
is a matrix of the form $(c^j_{rs} Q_j)_{r,s=1}^{m_j}$ for some
complex scalars $c^j_{rs}$. From these, a computation shows that the
entries of the matrix
$\tilde{u} = S diag(m_1 u_1, m_2 u_2, \cdots, m_k u_k) S^{-1}$
satisfy the defining relations for $A_u(Q)$,
where $u_j$ is the fundamental representation of $A_u(Q_j)$ 
($j=1, \cdots , k$).
Hence there is a surjection $\pi$ from $A_u(Q)$ to
$A_u(Q_1) \ast A_u(Q_2) \ast \cdots \ast A_u(Q_k)$
such that $\pi(u) = \tilde{u}$. That is $\pi(w_j) = u_j$
($j=1, \cdots , k$).

Again by Lemma 1.2 of \cite{W5} and the properties of free 
product Woronowicz $C^*$-algebras \cite{W1}, 
there is a surjection $\rho$ from
$A_u(Q_1) \ast A_u(Q_2) \ast \cdots \ast A_u(Q_k)$ to $A_u(Q)$
such that $\rho(u_j) = w_j$ ($j=1, \cdots , k$).
This is the inverse of $\pi$.
\QED
\vv
COROLLARY 1.
(1). {\em Let
$Q = diag(e^{i \theta_1} P_1, e^{i \theta_2} P_2, \cdots, e^{i \theta_k} P_k)$,
with positive matrices $P_j$ and
distinct angles $0 \leq \theta_j < 2 \pi$ ($j = 1, \cdots , k$, $k \geq 1$,
every normal matrix is unitarily equivalent to one such,
unique up to permutation of the indices $j$).
Then
\begin{center}
$A_u(Q) \cong A_u(P_1) \ast A_u(P_2) \ast \cdots \ast A_u(P_k).$
\end{center}}
(2). {\em Let $Q \in GL(2, {\mathbb C})$ be a non-normal matrix. Then
$A_u(Q) = C({\mathbb T})$.}

(3). {\em For $Q  \in GL(2, {\mathbb C})$, $A_u(Q)$ is either isomorphic 
to $C({\mathbb T})$, or $C({\mathbb T}) \ast C({\mathbb T})$, 
or $A_u(diag(1, q))$ with $0 < q \leq 1$. }
\vv
\pf
(1).  
Let $S$ and $E$ be as in the proof of Theorem 3. Since we do not have 
an explicit formula for the Haar measure of $A_u(Q)$, we must
determine the matrices $Q_j$ in Theorem 3 by other means.

We have an evident surjection 
$$\pi: A_u(Q) \rightarrow A_u(P_1) \ast A_u(P_2) \ast \cdots \ast A_u(P_k)$$ 
such that $( \pi (u_{ij}) ) = diag(u_1, \cdots , u_k) , $ 
where $u = (u_{ij})$ is the fundamental representation of $A_u(Q)$ 
and $u_j$ is the fundamental representation
of $A_u(P_j)$ ($j = 1,  \cdots , k$). That is, the free product quantum group
of the $A_u(P_j)$'s is a quantum subgroup of
$A_u(Q)$ (cf. \cite{W1} for the terminology).  Since the $u_j$'s
are mutually inequivalent representations (cf. Theorem 3.10 of \cite{W1}),
the multiplicities of the irreducible constituents of $u$ are all equal to
one and the matrix $E$ in the proof of Theorem 3 is of the form
\begin{center}
$E = diag(c_1 Q_1,  \cdots , c_l Q_l)$
\end{center}
for some $l \leq k$ ($c_1,  \cdots , c_l \in {\mathbb C}^*$, $Q_1,  \cdots , Q_l>0$).
Since the angles $\theta_j$ are distinct and $E$ is unitarily equivalent to
$Q$, we must have $l=k$,  $|c_j| Q_j = P_j$ and $c_j |c_j|^{-1} =
e^{ i \theta_j}$ after a possible permutation of the indices $j$
(Note that permutation of the indices $j$ does not change the quantum group
$B_u(T_1) \ast \cdots \ast B_u(T_k)$).
We conclude the proof by noting that $A_u(Q_j) = A_u(|c_j| Q_j)= A_u(P_j)$.

(2). Since $Q$ is not positive,
the fundamental representation $u$ of $A_u(Q)$ is reducible (cf.
3.1 of \cite{W14}). Since $Q$ is not normal, we deduce from
(1) and the proof of Theorem 3 that $u$ is equivalent to
a representation of the form $2 w_1$ (i.e. $m_1=2$), where $w_1$
an irreducible representation of dimension 1. 

(3). This  follows from (1) and (2) (cf. also Theorem 1). 
\QED
\vv
THEOREM 4.
{\em Let $Q \in GL(n, {\mathbb C})$. Then the fundamental representation $u$
of $B_u(Q)$ has a unitary isotypical decomposition of the form
\begin{center}
$u = S diag (m_1 w_1, m_1 \tilde{w}_1,
             m_2 w_2, m_2 \tilde{w}_2, \cdots, m_k w_k, m_k \tilde{w}_k,
             m'_1 w'_1, m'_2 w'_2, \cdots, m'_l w'_l) S^{-1}$
\end{center}
for some $S \in U(n)$, where the $w_i$'s are not self-conjugate,
$\tilde{w}_i = P_i^{1/2} \bar{w}_i P_i^{-1/2}$,
$P_i = h(w_i^t \bar{w}_i)$,  the $w'_j$'s are self-conjugate
($i=1, \cdots ,k$, $j=1,  \cdots , l$, $k, l \geq 0$). 
Let $Q_j$ be such that $w'_j = Q^*_j \overline{w'}_j {Q^*_j}^{-1}$.
Then the $Q_j \bar{Q}_j$'s are nonzero scalar matrices and
\begin{center}
$B_u(Q) \cong A_u(P_1) \ast A_u(P_2) \ast \cdots \ast A_u(P_k) \ast
B_u(Q_1) \ast B_u(Q_2) \ast \cdots \ast B_u(Q_l).$
\end{center}
}
\noindent
\pf
Note that the fundamental representation $u$ of $B_u(Q)$
is self-conjugate: $u = Q^* \bar{u} {Q^*}^{-1}$.
Hence, in the isotypical decomposition of $u$, if an irreducible
component is not self-conjugate, then its conjugate representation
also appears (with the same multiplicity as the former). By
Lemma 1.2 of \cite{W5}, we deduce that each $P_i^{1/2} \bar{w}_i P_i^{-1/2}$ is
a unitary representation.  Hence $u$ has a decomposition as stated in the
theorem. Since the $w'_j$'s are irreducible, $Q_j \bar{Q}_j$
are scalar matrices as in \cite{Banica1}.

Let $E = S^t Q S$.  Then the second set of 
relations for $B_u(Q)$ becomes 
$$ w^t E w E^{-1} = I_n = E w E^{-1} w^t, 
 \hspace{.51cm}
{\mbox i.e.} \hspace{.51cm}
w^t E w =  E, $$
where 
\begin{center}
$w =  diag ( v_1,  \tilde{v}_1,
             v_2,  \tilde{v}_2, \cdots, v_k,  \tilde{v}_k,
             v'_1, v'_2, \cdots, v'_l),$ 
\end{center}
$v_i = m_i w_i$, 
$\tilde{v}_i = m_i \tilde{w}_i$, 
$v'_j = m'_j w'_j$, $i = 1, \cdots , k, j = 1, \cdots, l$. 
Block decompose $E$ according to $w$,  we find that
$E$ is of the form
\begin{center}
$E = diag(
\left[
\begin{array}{cc}
  0 & X_1    \\
Y_1 & 0
\end{array}
\right],
\left[
\begin{array}{cc}
  0 & X_2    \\
Y_2 & 0
\end{array}
\right], \cdots ,
\left[
\begin{array}{cc}
  0 & X_k    \\
Y_k & 0
\end{array}
\right], Z_1, Z_2,  \cdots , Z_l )$,  
\end{center}
and the above set of relations takes the form 
$$ v_i^t X_i \tilde{v}_i = X_i, \; \; \; 
   \tilde{v}_i^t Y_i v_i = Y_i, \; \; \; 
   {v'_j}^t Z_j v_j = Z_j.$$
By the assumptions in Theorem 4, 
we have that (cf. proof of Theorem 3) 
$$(w_i^t)^{-1} = P_i \bar{w}_i P_i^{-1}, \; \; \;
  (\tilde{w}_i^t)^{-1} = 
   \bar{P}_i^{- \frac{1}{2}} w_i \bar{P}_i^{\frac{1}{2}}, \; \; \;
  ({w_j'}^t)^{-1} = Q_j w'_j Q_j^{-1} .$$ 
Hence we have 
\begin{center}
$X_i = (x_{rs}^i P_i^{1/2})_{r,s=1}^{m_i}, \hspace{5mm}
Y_i = (y_{rs}^i \bar{P}_i^{-1/2})_{r,s=1}^{m_i}, \hspace{5mm}
Z_j = (z_{rs}^j Q_j)_{r,s=1}^{m'_j},   \hspace{5mm} 
x_{rs}^i, y_{rs}^i, z_{rs}^j \in {\mathbb C}.$
\end{center}
\noindent
From these, a computation then shows that the entries of the matrix
\begin{center}
$\tilde{u} = S diag (m_1 u_1, m_1 \tilde{u}_1,
             m_2 u_2, m_2 \tilde{u}_2, \cdots, m_k u_k, m_k \tilde{u}_k,
             m'_1 u'_1, m'_2 u'_2, \cdots, m'_l u'_l) S^{-1}$
\end{center}
satisfy the defining relations for $B_u(Q)$,
where $u_i$ (resp. $u'_j$) is the fundamental representation
of $A_u(P_i)$ (resp. $B_u(Q_j)$) and
$\tilde{u}_i = P_i^{1/2} \bar{u}_i P_i^{-1/2}$
 ($i=1,  \cdots , k$, $j=1,  \cdots , l$).
Hence there is a surjection $\pi$ from $B_u(Q)$ onto
$ A_u(P_1) \ast A_u(P_2) \ast \cdots \ast A_u(P_k) \ast
B_u(Q_1) \ast B_u(Q_2) \ast \cdots \ast B_u(Q_l)$,
such that $\pi(u) = \tilde{u}$. That is
$\pi(w_i) = u_i$ and $\pi(w'_j) = u'_j$ ($i=1,  \cdots , k$, $j=1,  \cdots  ,l$).
As in the proof of Theorem 3, $\pi$ is an isomorphism.
\QED
\vv
COROLLARY 2.
(1). {\em Let
$Q = diag( T_1, T_2, \cdots, T_k)$ be a matrix such that
$T_j \bar{T}_j = \lambda_j I_{n_j}$, where the $\lambda_j$'s are distinct
non-zero real numbers
(the sizes $n_j$ need not be different), $j = 1,  \cdots , k$, $k \geq 1$.
Then \begin{center}
$B_u(Q) \cong B_u(T_1) \ast B_u(T_2) \ast \cdots \ast B_u(T_k).$
\end{center}}
\noindent
(2).
{\em Let
$Q =
\left[
\begin{array}{cc}
             0 & T    \\    
q \bar{T}^{-1} & 0
\end{array}
\right]$,
where $T \in GL(n, {\mathbb C})$ and $q$ is a complex but non-real number.
Then $B_u(Q)$ is isomorphic to $A_u(|T|^2)$ under the the map $\pi$ which
sends the entries of the fundamental representation
$u$ of $B_u(Q)$ to the entries of the matrix $diag(u_1, u_2)$, where
$u_1$ is the fundamental representation of $A_u(|T|^2)$ and
$u_2 = |T| \bar{u}_1 |T|^{-1}$, the unitary equivalent of $ \bar{u}_1$.
}
\vv
\pf 
(1).
Let $S$ and $E$ be as in the proof of Theorem 4.
We have an evident surjection from $B_u(Q)$ onto the free product
$ B_u(T_1) \ast B_u(T_2) \ast \cdots \ast B_u(T_k)$ sending the matrix
entries of the fundamental representation $u$ of $B_u(Q)$ to entries of
$diag(u_1, \cdots , u_k)$, where $u_j$ is the fundamental representation
of $B_u(T_j)$ ($j = 1,  \cdots , k$). Since the $u_j$'s
are mutually inequivalent self-conjugate representations and none of them
is conjugate to {\em another} (cf. Theorem 3.10 of \cite{W1}),
Theorem 4 implies that the pieces $A_u(P_i)$ do not appear in the
decomposition of $B_u(Q)$ and that the multiplicities $m'_j=1$.
Therefore the matrix $E$ in the proof of Theorem 4 has the form
\begin{center}
$E = diag( Z_1, Z_2,  \cdots , Z_l )$,
\end{center}
for some $l \leq k$, where
$ Z_j = z_j Q_j$, $ z_j \in {\mathbb C}^*$,
and the $Q_j \bar{Q}_j$'s are scalar matrices ($j=1, \cdots , l$).
Hence,
\begin{center}
$E \bar{E} = diag(c_1 I_{n'_1},  \cdots , c_l I_{n'_l})$,
\end{center}
for some $c_j \in {\mathbb C}^*$,
where $n'_j$ is the size of the matrix $Q_j$, $j=1,  \cdots , l$.
From the unitary equivalence
\begin{center}
$E \bar{E}
=S^t Q \bar{Q} \bar{S} =
S^t diag( \lambda_1 I_{n_1},  \cdots ,  \lambda_k I_{n_k} ) \bar{S}$
\end{center}
and the fact that the $\lambda_j$'s are distinct,
we must have $l = k$, $c_j = \lambda_j$ and
$n'_j = n_j$, up to a possible permutation of the indices $j$.
Now $E \bar{E} = S^t Q \bar{Q} \bar{S}$ takes the form
$$diag(c_1 I_{n_1},  \cdots , c_l I_{n_k}) =
S^t diag(c_1 I_{n_1},  \cdots , c_l I_{n_k}) \bar{S}.$$
Then from the assumption that the $c_j$'s are distinct we deduce again 
that $S$ is a block diagonal matrix $S = diag(S_1,  \cdots , S_k)$ with
$S_j \in U(n_j)$, $j = 1,  \cdots , k$.
Hence
\begin{center}
$E = S^t Q S = diag( S^t_1 T_1 S_1,  \cdots , S^t_k T_k S_k),$
\end{center}
and therefore $Z_j = z_j Q_j = S_j^t T_j S_j$, $j = 1,  \cdots , k$.
Hence by Theorem 2, $B_u(Q_j)$ is isomorphic to $B_u(T_j)$. 
The proof is finished by Theorem 4.

(2).
Let $T = U |T|$ be the polar decomposition of $T$. Then
\bgeqq
Q =
\left[
\begin{array}{cc}
             0 & T    \\    
q \bar{T}^{-1} & 0
\end{array}
\right]
=
\left[
\begin{array}{cc}
             U & 0    \\    
             0 & 1
\end{array}
\right]
\left[
\begin{array}{cc}
             0 & |T|    \\  
q \bar{|T|}^{-1} & 0
\end{array}
\right]
\left[
\begin{array}{cc}
             U^t & 0    \\    
              0  & 1
\end{array}
\right].
\ndeqq
By Theorem 2, we can assume $T > 0$ from now on. 
Then 
\bgeqq
(Q^*)^{-1} =
\left[
\begin{array}{cc}
                   0 & T^{-1}    \\    
\bar{q}^{-1} \bar{T} & 0
\end{array}
\right]. 
\ndeqq
Let $u = diag(u_1, u_2)$, where $u_1$ and $u_2$ are the unitaries as 
given in the statement of Corollary 2. Then,  a quick computation shows that  
\bgeqq
Q^* \bar{u} (Q^*)^{-1} = u,
\ndeqq
that is,
$
u^{-1} = u^* = Q^{-1} u^t Q. 
$
Hence we have a
surjection $\pi$ from $B_u(Q)$ onto $A_u(T^2)$ as
in the statement of Corollary 2.(2).
Then Theorem 3.10 of \cite{W1} and
Theorem 4 above implies that $B_u(Q)$ is isomorphic to
$A_u(P_1)$ for some $P_1>0$.
Now the rigidity of $A_u(P_1)$ (see Remark (1) in Sect. 1) implies that
$A_u(P_1)$, and therefore $B_u(Q)$, is isomorphic to $A_u(T^2)$.
\QED
\vv
{\em Concluding Remarks.}
(1). Using Theorem 1 (resp. Theorem 2) along with \cite{W1}, one can
show that the decomposition in Theorem 3 (resp. Theorem 4)
is unique in the evident sense.

(2).  Theorems 1, 2, 3, 4 solve Problem
1.1 of \cite{W7}.

(3). Using the same method as in Theorems 3 and 4, we
see that intersections of quantum groups of the form $A_u(Q)$
(resp. $B_u(Q)$) in the sense of \cite{W7} does not give rise to
non-trivial finite quantum groups. This solves Problem 2.4
of \cite{W7}.
\vv
{\bf Acknowledgment.}
Supported by NSF grant DMS-9627755. The author would  
like to thank the referee for pointing out a few obscurities in the 
original version of the paper.


\end{document}